\documentclass[12pt]{article}\usepackage{amsfonts}\usepackage{graphicx}
\usepackage{latexsym}\usepackage{amsbsy}

\begin{document}
\title{Embedding Bratteli-Vershik 
 systems in cellular automata\footnote{This research was partially supported by NSERC Canada.}}
\author{Marcus Pivato \ and \ Reem Yassawi \\
 Department of Mathematics, Trent University\thanks{1600 West Bank Drive,
Peterborough, Ontario,  K9J 7B8,  Canada.  Email:
 {\tt marcuspivato@trentu.ca} and {\tt ryassawi@trentu.ca}} }

\maketitle
 
\newcommand{\xinit}{\mathbf{x}_{\mathrm{init}}}

\newtheorem{theorem}{Theorem}
\newtheorem{proposition}[theorem]{Proposition}
\newtheorem{lemma}[theorem]{Lemma}
\newtheorem{corollary}[theorem]{Corollary}

\def\endproof{\hfill{\vrule height4pt width6pt depth2pt}

\vspace{0.5em}

}
\def\btau{{\boldsymbol{\tau}}}
\newsavebox{\clockbox}
\savebox{\clockbox}{\includegraphics[angle=-90,scale=0.35]{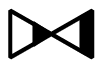}}
\newcommand{\clock}{\raisebox{0.7em}{\usebox{\clockbox}}}

\newcommand{\tilebox}[3][0.5]{\raisebox{-#2em}{\includegraphics[scale=#1]{#3.eps}}}

\begin{abstract}

  Many dynamical systems can be naturally represented as {\em
Bratteli-Vershik} (or {\em adic}) systems, which provide an appealing
combinatorial description of their dynamics.  If an adic system $X$
satisfies two technical conditions ({\em focus} and {\em bounded
width}) then we show how to represent $X$ using a two-dimensional
subshift of finite type $Y$; each `row' in a $Y$-admissible
configuration corresponds to an infinite path in the Bratteli diagram
of $X$, and the vertical shift on $Y$ corresponds to the `successor'
map of $X$.  Any $Y$-admissible configuration can then be recoded as
the spacetime diagram of a one-dimensional cellular automaton $\Phi$;
in this way $X$ is {\em embedded} in $\Phi$ (i.e. $X$ is conjugate to
a subsystem of $\Phi$).  With this technique, we can embed many
odometers, Toeplitz systems, and constant-length substitution systems
in one-dimensional cellular automata.

\end{abstract}

\section{Introduction}

  The Bratteli-Vershik (`adic') transformations are a large and
diverse class of symbolic dynamical systems which are of interest
partly because many other dynamical systems can be naturally
represented as adic transformations, thereby yielding an appealing
combinatorial description of their dynamics.  Cellular automata (CA)
are another class of symbolic dynamical system with versatile
representation capabilities.  Recent work shows that certain simple adic
systems arise naturally as subsystems of CA.  For example, if
$\Phi:{\cal A}^{\mathbb Z}\rightarrow{\cal A}^{\mathbb Z}$ is a
non-injective, right-sided, left-permutative CA, then many
orbit-closures of $\Phi$ are conjugate to odometers \cite{cpy}; a more
precise characterization of some of these odometers appears in
\cite{cy}.  Similarly, in the {\em Game of Life} CA, the orbit closure
of a certain configuration can be represented as a substitution system
\cite{gm}.  The main result of this article is the following:

   

\begin{theorem}
\label{mainthm}
Let $(X_{\cal B}, V_{\cal B})$ be an adic system, where
${\cal B}$ is a  properly ordered Bratteli diagram which is focused, and has bounded width. 
Then there exist 
cellular automata which embed $(X_{\cal B}, V_{\cal B})$.
\end{theorem}

 For definitions of terms in Theorem
\ref{mainthm}, and notation, see Section \ref{introsection}. Theorem \ref{mainthm} follows directly from Theorem \ref{thisisit} (in Section \ref{proof}).
In Section \ref{examples} we give examples
of families of dynamical systems whose adic representations 
satisfy the conditions of Theorem \ref{mainthm}. This includes all primitive substitutions, and many Toeplitz
systems (Corollaries \ref{cor1} and \ref{toeplitzcorollary}).  This technique also reproves the
embedding of certain odometers in a CA, a result in \cite{cpy};
however the family of CA that we identify here are disjoint from
those found in \cite{cpy}.

\section{Preliminaries \label{introsection}}
\subsection{Notation} 
 Let ${\cal A}$ be a finite alphabet; elements in ${\cal A}$ will be
denoted by $a,\, b\,, c\,$ etc.  A {\em word} ${\bf a}$ from ${\cal
A}$ is a finite concatenation of elements from ${\cal A}$.  Let ${\cal
A}^{k}$, ${\cal A}^{\leq k}$, and ${\cal A}^{+}$ denote the set of all
words of length $k$, length at most $k$, and any finite length
respectively. We include the empty word, and use boldface to denote
words.  If ${\bf a}=a_{0}a_{1}\ldots a_{j}$ and ${\bf
a^{*}}=a_{0}^{*}a_{1}^{*}\ldots a_{k}^{*}$, then ${\bf aa^{*}}:=
a_{0}a_{1} \ldots a_{j}a_{0}^{*}a_{1}^{*}\ldots a_{k}^{*}$. If ${\bf
a}= a_{0}a_{1}\ldots a_{j}$, define $|{\bf a}|:= j$.  The space of all
bi-infinite sequences from ${\cal A}$ is written as ${\cal A}^{\mathbb
Z}$. Elements of the latter are written ${\bf a} =\ldots a_{-1}\cdot
a_{0}a_{1} \ldots $.  If ${\mathbb N} := \{0,1,2,\ldots\},$ and
${\mathbb N}^{+} := \{1,2,\ldots\},$ elements of ${\cal A}^{\mathbb
N}$ and ${\cal A}^{\mathbb N^{+}}$ are defined analagously.  The
notation $\{x_{n}\}_{n=\infty}^{k}$ is used to denote the
left-infinite sequence $\ldots x_{k+2}\,x_{k+1}\, x_{k}$. If ${\bf x}
= \{x_{n}\}_{n=\infty}^{1}$ and ${\bf y} = \{y_{n}\}_{n=0}^{\infty}$,
then by ${\bf x}\cdot {\bf y}$ we mean the bi-infinite sequence
$\ldots x_{2}\,x_{1}\cdot y_{0}\, y_{1} \ldots$.  If ${\mathbb L}=
{\mathbb Z}, \,\,\, {\mathbb N},$ or ${\mathbb N}^{+}$, then ${\cal
A}^{\mathbb L}$ is a {\em Cantor space}: a  zero-dimensional compact metric space, when ${\cal A}$ is endowed
with the discrete topology and ${\cal A}^{\mathbb L}$ with the product
topology. If $b\,\in {\cal A}$, let $[b]_{j}:=\{{\bf a}:
\,\,a_{j}=b\}$; these (clopen) sets form a countable basis for the
topology on ${\cal A}^{\mathbb L}$.  This topology is also generated
by the Hamming metric.  The (left) {\em shift map} $\sigma: {\cal
A}^{\mathbb L} \rightarrow {\cal A}^{\mathbb L}$ is the map defined as
$(\sigma(x))_{n} = x_{n+1}$.  X is a {\em subshift of $({\cal
A}^{\mathbb L},\sigma)$} if it is a closed $\sigma$-invariant subset
of ${\cal A}^{\mathbb L}$.
A (one dimensional) {\em cellular automaton} is a continuous,
$\sigma$-commuting map $\Phi: {\cal A}^{\mathbb Z} \rightarrow {\cal
A}^{\mathbb Z}$.  The Curtis-Hedlund-Lyndon theorem \cite{h} states
that every CA is given by a {\em local rule} $\phi :{\cal A}^{l+r+1
}\rightarrow {\cal A}$ for some $l\geq 0$ (the {\em left radius of
$\phi$}) and $r\geq0$ (the {\em right radius of $\phi$)}, where for
all $x\,
\in {\cal A}^{\mathbb Z}$, and all $i \, \in {\mathbb Z}$,
\[
[\Phi(x)]_{i} = \phi(x_{i-l}, x_{i-l+1}, \ldots , x_{i+r}) \, . \]

If $X\subset {\cal A}^{\mathbb Z}$ is closed and $\Phi(X)\subseteq X$,
then $(X,\Phi)$ is a {\em subsystem} of $({\cal A}^{\mathbb
Z},\Phi)$. If the dynamical system $(Y,T)$ is topologically conjugate
to the subsystem $(X,\Phi)$ of $({\cal A}^{\mathbb
Z},\Phi)$ (i.e. if  there exists a homeomorphism $f: Y\rightarrow X$ with 
$f\circ T= \Phi\circ f$), we say
that $({\cal A}^{\mathbb Z},\Phi)$ {\em embeds} $(Y,T)$.

\subsection{Bratteli Diagrams}
A {\em Bratteli diagram} ${\cal B}=({\cal V},{\cal E})$ is an infinite
directed graph with {\em vertex set} ${\cal V} = \bigsqcup_{n=
0}^{\infty} {\cal V}_{n}$ and {\em edge set} ${\cal E} =
\bigsqcup_{n=1}^{\infty}{\cal E}_{n}$, where all ${\cal V}_{n}$'s and
${\cal E}_{n}$'s are finite, ${\cal V}_{0}= \{\clock\}$, and, if $x$
is an edge in ${\cal E}_{n}$, the source $s(x)$ of $x$ lies in ${\cal
V}_{n}$ and the range $r(x)$ of $x$ lies in ${\cal V}_{n-1}$. We
remark that in all references to Bratteli diagrams that we have
consulted, edges move in the opposite direction: they have source in ${\cal V}_n$ and range ${\cal
V}_{n+1}$; however for the purposes of our results in Section \ref{proof}, we find our representation more visually intuitive.  We will use $x, y\ldots$ when referring to edges,
and $a,b, \ldots$ when referring to vertices, although edges in ${\cal
E}_{1}$ will be given special labels $\clock_{i}$, where $i$ varies.
A finite set of edges $\{x_{n+k}\}_{k=1}^{K}$, with
$s(x_{n+k})=r(x_{n+k+1})$ for $1\leq k\leq K-1$, is called a {\em path
} from $s(x_{n+K})$ to $r(x_{n+1})$. Similarly an {\em infinite path}
in ${\cal B}$ is a sequence ${\bf x}=\{x_{n}\}_{n=\infty}^{1}$, with
$x_{n}\,\in {\cal E}_{n}$ for $n\geq 1$, and $s(x_{n})=r(x_{n+1})$ for
$n\geq 1$.  We will often write
\begin{equation}
\label{visual}
\ldots a_{n+1}\stackrel{x_{n+1}}{\rightarrow}a_{n}
\stackrel{x_{n}}{\rightarrow}\ldots \stackrel{x_{2}}{\rightarrow} a_{1}
\stackrel{x_{1}}{\rightarrow} \clock
\end{equation}
where $s(x_{n})=a_{n}$.  The set of all infinite paths in ${\cal B}$
will be denoted $X_{\cal B}$ (a subset of $ \Pi_{n\geq 1} {\cal E}_{n}
$), and  $X_{\cal B}$ is endowed with the topology induced from the product
topology on $\Pi_{n\geq 1} {\cal E}_{n}$. Thus $X_{\cal B}$ is a
compact metric space.

If ${\bf x}= \{x_{n}\}_{n=\infty}^{1}$ and ${\bf
x'}=\{x_{n}'\}_{n=\infty}^{1}$ are two elements in $X_{\cal B}$, we
write ${\bf x} \sim {\bf x'}$ if the tails of ${\bf x}$ and ${\bf x'}$
are equal. It follows that $\sim$ is an equivalence relation, and primitivity of $\tau$
implies that each equivalence class of $\sim$ is dense.  We will mostly write $({\bf a},{\bf x}) =
\{(a_{n},x_{n})\}_{n\geq 1}$, where $s(x_{n})=a_{n}$, when referring to
an element ${\bf x}$ in $X_{\cal B}$.

Two Bratteli diagrams ${\cal B}=({\cal V}, {\cal E})$ and ${\cal B'}=({\cal V'}, {\cal E'})$ are {\em isomorphic} if 
there exists a pair of bijections $f_{{\cal V}}: {\cal V}\rightarrow {\cal V'}$ and 
$f_{\cal E}:{\cal E}\rightarrow {\cal E'}$ satisfying $f_{\cal V}(a) \, \in {\cal V}_{n}'$ if $a \, \in {\cal V}_{n}$,
and $s(f_{\cal E}(x)) = f_{\cal V}(s(x))$,  $r(f_{\cal E}(x)) = f_{\cal V}(r(x))$ whenever $x\, \in {\cal E}$.

Let $\{n_{k}\}_{k=0}^{\infty}$ be a sequence of increasing integers
with $n_{0}=0$. Then ${\cal B}' =({\cal V}',{\cal E}')$ is a {\em
telescoping} of ${\cal B}= ({\cal V},{\cal E})$ if ${\cal
V}_{k}'=V_{n_k}$ (with the vertex $v \, \in {\cal V}_{n_k}$ labelled
as $v' \, \in {\cal V}_{k}'$), and the number of edges from
$v_{k+1}'\, \in {\cal V}_{k+1}'$ to $v_{k}'\, \in {\cal V}_{k}'$ is
the number of paths from $v_{k+1}\, \in {\cal V}_{n_{k+1}}$ to
$v_{k}\, \in {\cal V}_{n_{k}}$.  We consider two Bratteli diagrams
${\cal B}$ and  ${\cal B}'$ {\em equivalent} if ${\cal B'}$ can be
obtained from ${\cal B}$ by isomorphism and telescoping. Thus when we
talk about a Bratteli diagram we are talking about an equivalence
class of diagrams.
 We say that ${\cal B}$ is
{\em simple} if there exists a telescoping ${\cal B'}=({\cal V'},{\cal
E'})$ of ${\cal B}$ so that, for any $a \, \in {\cal V}'_{n+1}$ and
$b\, \in {\cal V}'_{n}$, there is at least one edge from $a$ to $b$.
If $X_{\cal B}$ is simple, then $X_{\cal B}$ has no isolated points,
making it a Cantor space.

\subsubsection{Ordering  $X_{\cal B}$}
For each $a \, \in {\cal V}_{n}$, let ${\cal E}_{n}(a)= \{x \, \in
{\cal E}_{n}: s(x)=a\}.$ Say ${\cal B}$ is {\em ordered } if for each
$n \geq 1$ and $a \, \in {\cal V}_{n}$, there is a linear order $\geq$
on ${\cal E}_{n} (a)$; elements of ${\cal E}_n$ will then be labelled $0,1,
\ldots$ according to their order. If $a\, \in {\cal V}\backslash
\{\clock\}$, define $|a|:= |{\cal E}_{n}(a)|-1$, so that ${\cal E}_{n}(a) = \{0,1, \ldots |a|\}$.

The linear order on
edges in each ${\cal E}_{n}(a)$ induces a partial ordering on paths from
${\cal V}_{n}$ to ${\cal V}_{m}$: the two paths ${\bf x}=
a_{n}\stackrel{x_{n}}{\rightarrow}a_{n-1}
\stackrel{x_{n-1}}{\rightarrow}\ldots
\stackrel{x_{m+1}}{\rightarrow}a_{m}$ and ${\bf
x'}=a_{n}'\stackrel{x_{n}'}{\rightarrow}a_{n-1}'
\stackrel{x_{n-1}'}{\rightarrow}\ldots
\stackrel{x_{m+1}'}{\rightarrow}a_{m}'$ from ${\cal E}_{n}$ to ${\cal
E}_{m}$ are comparable with ${\bf x}<{\bf x'}$ if $a_{n}=a_{n}'$ and
if there is some $k\, \in [n,m]$ with $x_{k}<x_{k}'$ and
$x_{j}=x_{j}'$ for $k+1\leq j\leq n$.

Finally, two elements ${\bf x},\, {\bf x'}\,\in X_{\cal B}$ are {\em comparable} with
${\bf x}<{\bf x'}$ if there is a $k$ such that $x_{n}=x_{n}'$
for all $n>k$,  and $ x_{k}<x_{k}'$. Thus each
equivalence class for $\sim$ is ordered.
There is the obvious notion of
ordered isomorphism of two ordered Bratteli diagrams ${\cal B},$
${\cal B'}$: the isomorphism between ${\cal B}$ and ${\cal B'}$ also
has to satisfy $f_{\cal E}(x) \leq f_{\cal E}(y)$ if $x\leq y$. If
${\cal B}'$ is a telescoping of the ordered Bratteli diagram ${\cal
B}$, then the order induced on ${\cal B}'$ from the order on ${\cal B}$
makes ${\cal B}'$ an ordered Bratteli diagram.  We say that the ordered
Bratteli diagrams ${\cal B},$ ${\cal B'}$ are {\em equivalent} if
${\cal B'}$ is the image of ${\cal B}$ by telescoping and order
isomorphism.

An infinite path is {\em maximal} ({\em minimal}) if all the edges
making up the path are maximal (minimal).  If ${\bf x}=
\{x_{n}\}_{n=\infty}^{1}$ is not maximal, let $k$ be the smallest
integer such that $x_{k}$ is not a maximal edge, and let $y_{k}$ be
the successor of $x_{k}$. Then the {\em successor} of $x$ is the path
$\ldots x_{k+2}\,x_{k+1}\, y_{k}\, 0 \ldots 0$. Similarly, every
non-minimal path has a {\em predecessor}.  Let $X_{\min}$, $(X_{\max})
\subset X_{\cal B}$ be defined as the set of minimal (maximal)
elements of $X_{\cal B}$. By compactness, these sets are non empty.
Let $V_{\cal B}:X_{\cal B}\backslash X_{max}\rightarrow X_{\cal
B}\backslash X_{min}$ be the successor map.
Simple ordered Bratteli diagrams which have a unique minimal and
maximal element (called ${\bf x}_{\min}$ and ${\bf x}_{\max}$ respectively) are
called {\em properly ordered}.  If ${\cal B} $ is properly ordered,
then $V_{\cal B}$ can be extended to a homeomorphism on $X_{\cal B}$
by setting $V_{\cal B} ({\bf x}_{\max})= {\bf x}_{\min}$.  We call $(X_{\cal
B},V_{\cal B})$ the {\em Bratteli-Vershik} or {\em adic} system associated 
with $\tau$. Note that $(X_{\cal
B},V_{\cal B})$ is a minimal system, since $V_{\cal B}$
orbits are equivalence  classes for $\sim$. Let us say that two Cantor
systems $(X_{i},T_{i},{\bf x}_{i})$, $i=1,2$ are {\em pointedly
isomorphic} if there exists a homeomorphism $f:X_{1}\rightarrow X_{2}$
with $f\circ T_{1}= T_{2}\circ f$ and $f(x_{1})=x_{2}$.  The Bratteli-
Vershik system associated to an equivalence class of properly ordered
Bratteli diagrams is well defined up to pointed isomorphism:

\begin{theorem}
\label{putnam}
Let ${\cal B}$ and ${\cal B}'$ be properly ordered Bratteli
diagrams. Then ${\cal B}$ is equivalent to ${\cal B}'$ if and only if
$(X_{\cal B}, V_{\cal B}, {\bf x_{\min}})$ is pointedly isomorphic to
$(X_{\cal B'}, V_{\cal B'}, {\bf x_{\min}})$.
\hfill{\rm\cite[\S4]{hps}}
\end{theorem}

We say that ${\cal B}$ has {\em bounded width } if there exists a
constant $K$ such that  $|{\cal V}_{n}|\leq K$ and $|{\cal E}_{n}|\leq K$ for each $n\geq 1$. If ${\cal B}$ is ordered, we say that ${\cal
B}$ is {\em focused } if, for each $n\geq 1$, all minimal edges in
${\cal E}_{n}$ have the same range.  For example, any {\em proper substitution} (see Section \ref{examples}) has a focused representation. If ${\cal B}$ is focused, we use
$ \underline{a}$ to denote the range of any minimal edge in
${\cal E}_{n}$. Note that if ${\cal B}$ is focused, then its unique  minimal element is $( {\underline{\bf a}},{\bf 0})$.

\section{The spacetime diagrams of $(X_{\cal B}, V_{\cal B})$ and its associated subshift of finite type \label{proof}}
In this section, $\cal B$ is a properly ordered Bratteli diagram,
which has bounded width K. In this case vertices in ${\cal V}$ can be
labelled from a finite alphabet ${\cal A}$. For each $n \geq 1$ there
is a function ${\boldsymbol{\tau}}^{n}: {\cal A} \rightarrow {\cal
A}^{\leq K}$ such that ${\btau}^{n}$ completely describes ${\cal
E}_{n}$ and its ordering. In other words, if we  write ${\btau}^{n}(a) = \tau^{n}_{0}(a)
\tau^{n}_{2}(a) \ldots \tau^{n}_{ l^n_a}(a)$ where $l_a^n =
|\btau^{n}(a)|$, 
for each $a\, \in {\cal A}$ and $n\geq 1$, then
there is an edge from $a \, \in {\cal V}_{n}$ to $b$ in ${\cal
V}_{n-1}$, and it is labelled $k$, if and only 
$\tau^{n}_{k}(a) = b$.
  If ${\btau}^{n}_{a}$ is the empty word,
this is taken to mean that $a$ does not appear as a vertex in ${\cal
V}_{n}$. Implicit in this description is some arbitrary but
pre-assigned labelling of vertices in ${\cal V}$ from ${\cal A}$. We
will also assume that if ${\cal B}$ is focused, then its labelling
reflects this.  We shall use `stars' such as $*, \star, \ldots$ to
indicate $l^n_a$, so that ``$\tau^{n}_{*}(a)$'' means the last letter of
${\btau}^{n}(a)$, provided that there is no ambiguity.

Define ${\cal A}_{\cal B}$ to be the alphabet
\begin{eqnarray*}
{\cal A}_{\cal B}& := &\{(a,{\btau}^{n}, x): a \, \in {\cal A} , n \, \in 
{\mathbb N},
 x\, \in [0, l^n_a] 
\}
\\
& \cup & \{(a,{\btau}^{1}, \clock_{x}): a \in {\cal A},  x\, \in [0,l^1_a]
\}\,\cup\{\clock\}\, .
\end{eqnarray*}

An element $(\mathbf{a,x}) \, \in X_{\cal B}$ can also be seen as an
 element of ${\cal A}_{\cal B}^{\mathbb N}$ by writing $(\mathbf{a,x}) =
 ({\bf a},\{ {\btau}^{n}\}_{n\geq 1}, {\bf x})$. Conversely, when we write ``$(
 \mathbf{a, x}) \, \in {\cal A}_{\cal B}^{\mathbb N}$'', we mean $( \mathbf{a,
 x})=(\mathbf{ a},\{ {\mathbf{\btau}}^{n}\}_{n\geq 1}, \mathbf{x})$.  
 Let $\tilde{Y}\subset {\cal
 A}_{\cal B}^{{\mathbb Z}\times {\mathbb N}^{+}}$ be the set of all
 {\em space-time diagrams} for $(X_{\cal B},V_{\cal B})$: elements
 ${\bf \tilde{y}} =\{{\bf \tilde{y}}_{m}\}_{m\in {\mathbb Z}} \, \in
 \tilde{Y}$ (with the row ${\bf \tilde{y}}_{m}=
 \{\tilde{y}^{m}_{n}\}_{n=\infty}^{1}$) are such that ${\bf
 \tilde{y}}_{0} ={\bf (a, x)}$ for some
 ${\bf(a,x)} \, \in X_{\cal B}$ and ${\bf \tilde{y}}_{m} =\{ V_{\cal
 B}^{m}{\bf (a,x)}\}$ for each $m\, \in {\mathbb Z}$. Let
 $\clock^{\infty} =  \clock,\clock,\ldots $; We extend
 $\tilde{Y}$ to $Y\subset {{\cal A}_{\cal B}}^{\mathbb Z \times
 \mathbb Z}$ by letting
\[
Y:= \{ {\bf y} = \{{\bf y}_{m}\}_{m\in {\mathbb Z}}: 
{\bf y}= {\bf \tilde{y}}\cdot {\clock}^{\infty} \mbox{ for some } {\bf \tilde{y}}
 \in \tilde{Y} \}\, .
\] 
There is a natural identification of elements in $Y$ with elements in $\tilde{Y}$.
Note that we use positive integers to denote column locations in ${\bf y}$, unless there is a possibility of confusion.

A (two-dimensional) {\em subshift} is a closed, shift-invariant
subset $\mathcal{S}\subset\mathcal{A}^{\mathbb{Z}^2}$.  If $\mathcal{S}\subset\mathcal{A}^{\mathbb{Z}^2}$ is a subshift, then for
any finite $\mathbb{K}\subset\mathbb{Z}^2$, let
 $\mathcal{S}_\mathbb{K}:=\{\mathbf{s}_\mathbb{K} \ ; \ \mathbf{s}\in\mathcal{S}\}$ be
the set of all {\em $\mathcal{S}$-admissible $\mathbb{K}$-blocks}.  If $\mathsf{z}\in\mathbb{Z}^2$,
then $\mathcal{S}_\mathbb{K}=\mathcal{S}_{\mathbb{K}+\mathsf{z}}$ (because $\mathcal{S}$ is shift-invariant).  We
say $\mathcal{S}$ is a {\em subshift of finite type} (SFT) if there is some
finite {\em neighbourhood} $\mathbb{K}\subset\mathbb{Z}^2$ such that
$\mathcal{S}=\{\mathbf{a}\in\mathcal{A}^{\mathbb{Z}^2} \ ; \ \mathbf{a}_{\mathbb{K}+\mathsf{z}}\in\mathcal{S}_\mathbb{K}, \ \forall
\mathsf{z}\in\mathbb{Z}^2\}$.  For example, if $\Phi:\mathcal{A}^\mathbb{Z}
\rightarrow\mathcal{A}^\mathbb{Z}$ is a one-dimensional cellular automaton, then the
set of spacetime diagrams for $\Phi$ is an SFT in $\mathcal{A}^{\mathbb{Z}^2}$.

\begin{figure}
\centerline{\includegraphics[angle=90,scale=0.5]{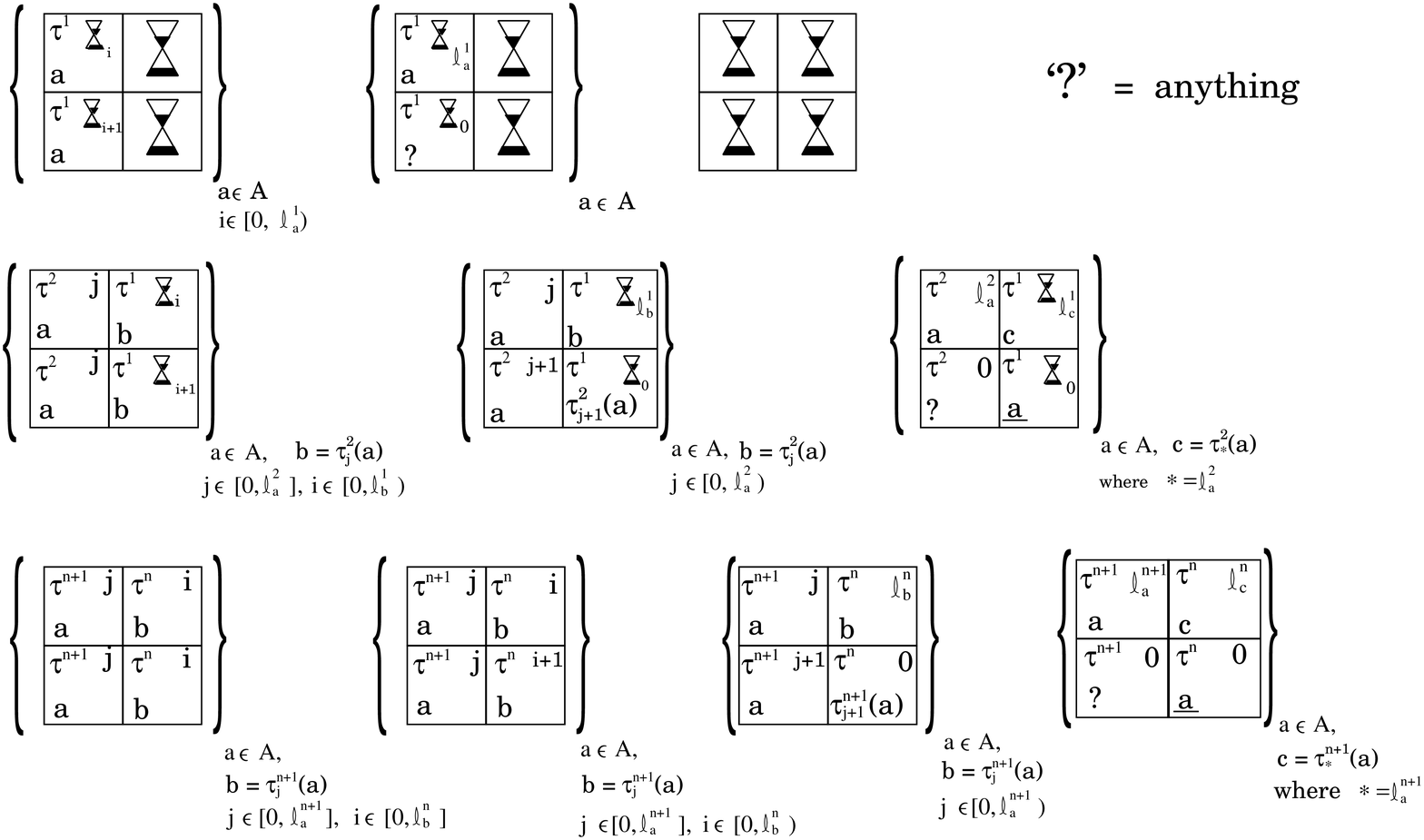}}
\caption{\footnotesize Admissible $\{0,1\}^2$-blocks defining the subshift
of finite type $\mathcal{S}$. \label{allowedtiles}}
\end{figure}

Let  $\mathbb{K}=\{0,1\}^2$, and let $\mathcal{S}$ be the SFT defined using the $\mathbb{K}$-blocks in Figure \ref{allowedtiles}.
 The following lemma is straightforward.
\begin{lemma}
\label{softlemma}
$Y$ is the  subset \
$\displaystyle \bigcup_{\mbox{$\scriptstyle a\,\in {\cal A}$}\atop\mbox{$\scriptstyle 0\leq i \leq l^{1}_a$}} \{ {\bf y} \in \mathcal{S} \ ; \  y_{0}^{0}= \clock, \
y_{-1}^{0} = (a,\clock_{i}) \}$
\ of $\mathcal{S}$.
\endproof\end{lemma}

We will sometimes use the symbols 
\tilebox{0.7}{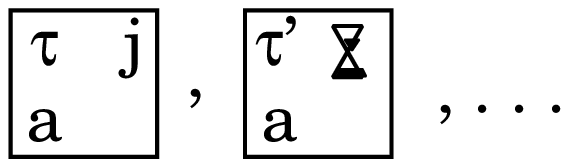}
or 
\tilebox{0.45}{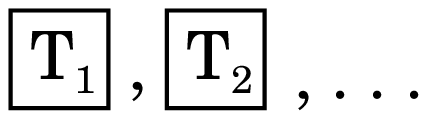}
to denote elements of ${\cal A}_{\cal B}$. For example, if ${\cal B}$ is the 
Bratteli diagram from Figure \ref{Brattelidiagram}, then
\[
{\cal A}_{\cal B}\quad=\quad
\tilebox{0.7}{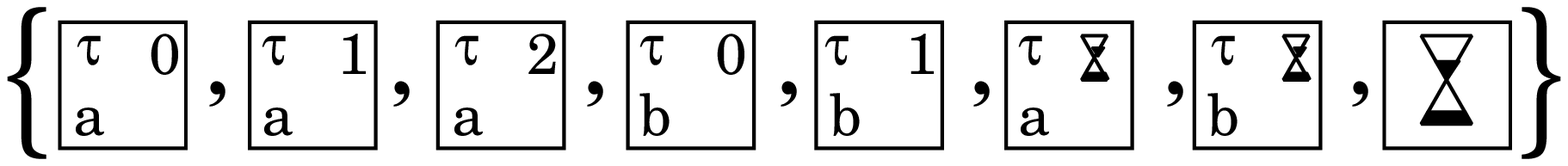},
\]
where the dropped index in $\clock$ indicates that there is a unique edge from each vertex in ${\cal V}_1$ to $\clock$.
A (possibly infinite) concatenation of letters from ${\cal A}_{\cal B}$ is a {\em
fragment} of a spacetime diagram if it appears in some ${\bf y}\, \in
Y$.

\begin{lemma}
\label{calemma}
If $\cal B$ is focused and properly ordered, and 
\tilebox[0.45]{1}{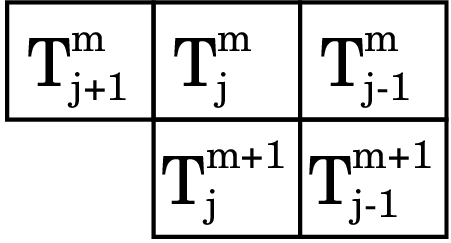}
is a
fragment of  ${\bf y} \, \in Y$, then 
\tilebox[0.45]{1}{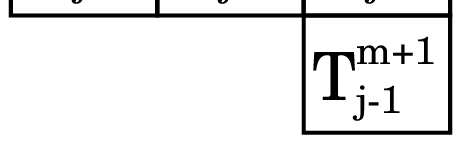}
determines
\tilebox[0.45]{0.5}{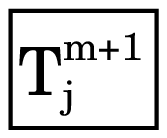}
.
\end{lemma}

{\bf Proof:} We prove the case where all tiles in
\tilebox[0.45]{1}{fig4}
are from 
$\{\tilebox[0.45]{0.5}{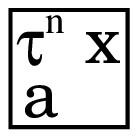}: a \, \in {\cal A} , n\geq 1 ,
 x\, \in [0, l^n_a]  \}$; other 
 cases are similar. If $T_{j-1}^{m}=T_{j-1}^{m+1}$, then
$T_{j}^{m+1}=T_{j}^{m}$, as is the case if $T_{j-1}^{m} = \tilebox[0.45]{0.5}{newfig1}$ with
$x<l^n_a$. If 
\[
\tilebox[0.65]{1.6}{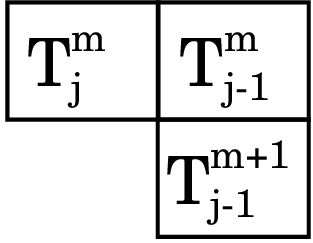}
\quad=\quad
\tilebox[0.55]{1.9}{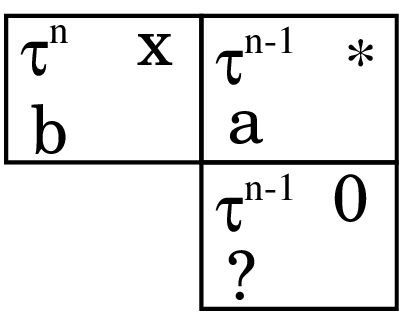}\]
with $x < l^n_b$, 
then $T_{j}^{m+1} = \tilebox[0.45]{0.5}{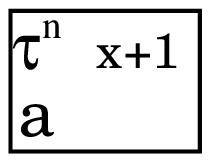}$. If 
\[
\tilebox[0.65]{1.6}{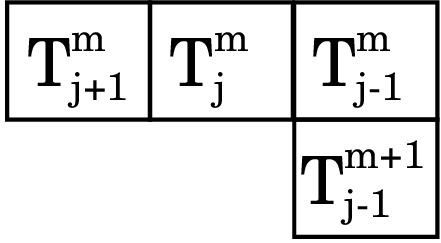}
\quad=\quad 
\tilebox[0.55]{1.9}{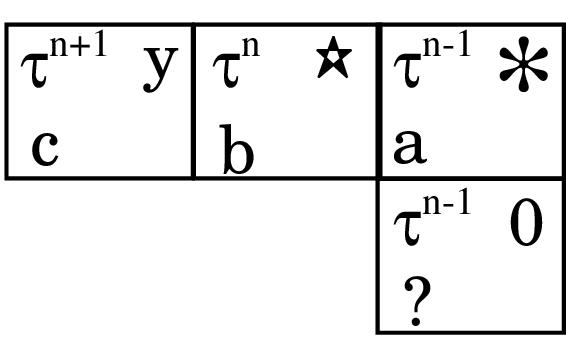}
\]
with $y<l^{n+1}_c$, then $T_{j}^{m+1} =\tilebox[0.45]{0.5}{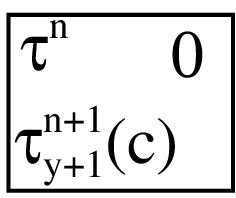}$. Finally if
\[
\tilebox[0.65]{1.6}{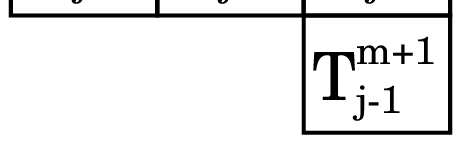}
\quad=\quad 
\tilebox[0.55]{1.9}{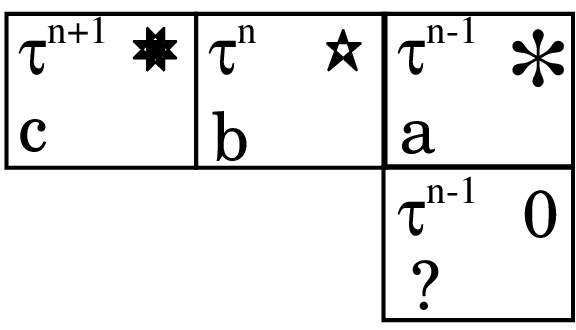}\,,\]

 then since $\tau$ is focused, we must have  $T_{j}^{m+1}=\tilebox[0.45]{0.5}{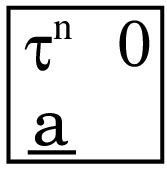}$. 
\endproof


Let  $W_{4}$ be the set of all fragments 
\tilebox[0.45]{1}{fig5}
appearing in $Y$, where $m,\, j,\, \in {\mathbb Z}$. Let 
$f:W_{4}\rightarrow :{\cal A}_{\cal B}$ be the
function such that if 
\tilebox[0.45]{1}{fig4}
is the unique extension of
\tilebox[0.45]{1}{fig5}, then 
$f\left(\mbox{\tilebox[0.45]{1}{fig5}} \right)= 
\mbox{\tilebox[0.45]{0.5}{fig6}}\,$.

 This also means that any (diagonal) ray
\[{\bf R}:=\{(\tilebox[0.45]{0.5}{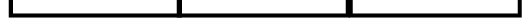}    ) : k\, \in {\mathbb N}\}\]  which is a fragment in $Y$
determines the horizontal ${\bf T}_{0}:=\{T_{k}^{0}: k \, \in {\mathbb
N}^{+}\}$. Let ${\cal R}$ be the collection of such rays ${\bf R}$,
and let $F: {\cal R}\rightarrow {\cal A}_{\cal B}^{\mathbb N}$ be the
function which maps these rays in ${\cal R}$ to their horizontal image
${\bf T}_{0}$. Call any 3-tuple $\tilebox[0.45]{0.5}{fig13}$ in a ray a
{\em step}. Note that any ray in ${\cal R}$ can be extended in a
unique way to a ``two-sided'' ray fragment from $Y$: simply attach
infinitely many steps of the form $\tilebox[0.45]{0.5}{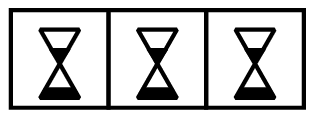}$.

Let ${\cal A}_{\Phi}:={\cal A}_{\cal B}^{3}$ and $W_{3}\subset {\cal
A}_{\Phi}$ be the set of all horizontal 3-tuple fragments from $Y$. We
will use $S_{0}, S_{1} \ldots$ to denote letters from ${\cal
A}_{\Phi}$. Define a left radius  one, right ra
3-++
dius one CA $\Phi:
W_{3}^{\mathbb Z}\rightarrow W_{3}^{\mathbb Z}$ with local rule
$\phi:W_{3}^{3}\rightarrow W_{3}$ defined as
\[
\tilebox{0.8}{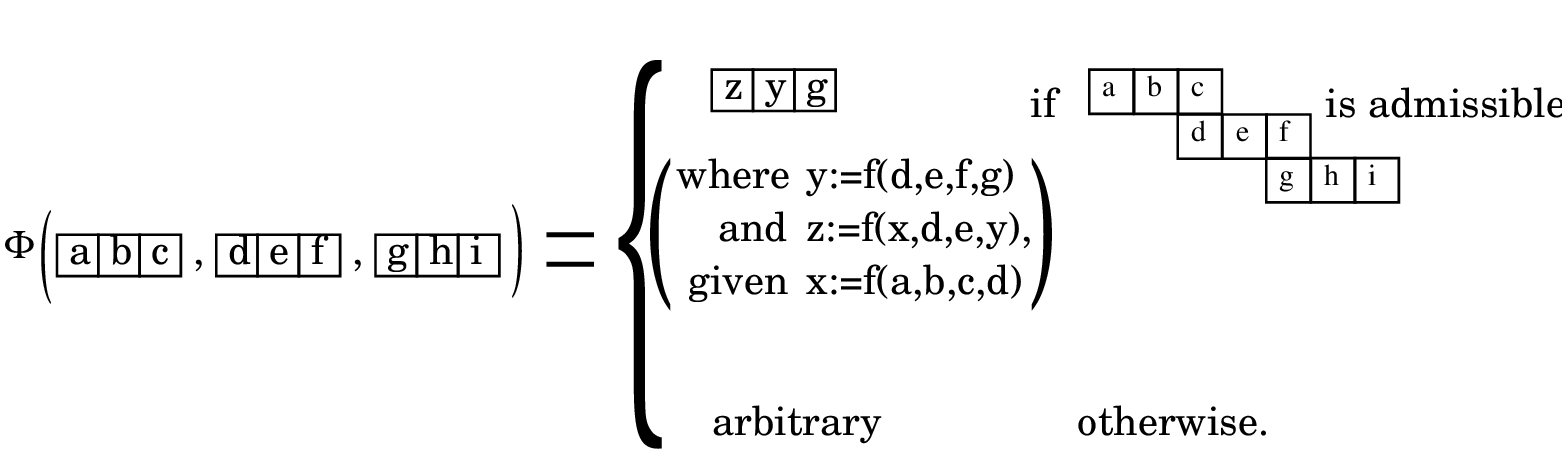}
\]
Let the maximal element of $X_{\cal B}$ be $(\{m_{n}\}_{n\geq 1}, 
\{ \star_{n}\}_{n\geq1})$, and  
define  
\begin{eqnarray*}
\lefteqn{\xinit \quad := }
\\ 
&& \includegraphics[angle=-90,scale=0.43]{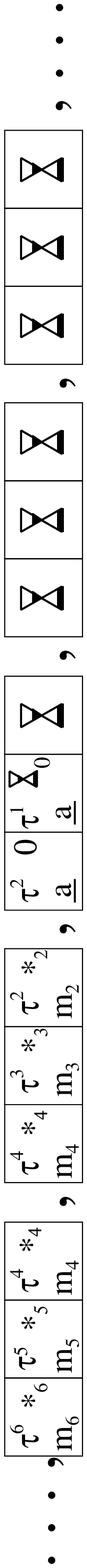}
\end{eqnarray*}
Set  $\Omega=\overline{\{\Phi^{n}((\xinit)\}_{n\geq 0}}$.

\begin{figure}
\centerline{\includegraphics[scale=0.6]{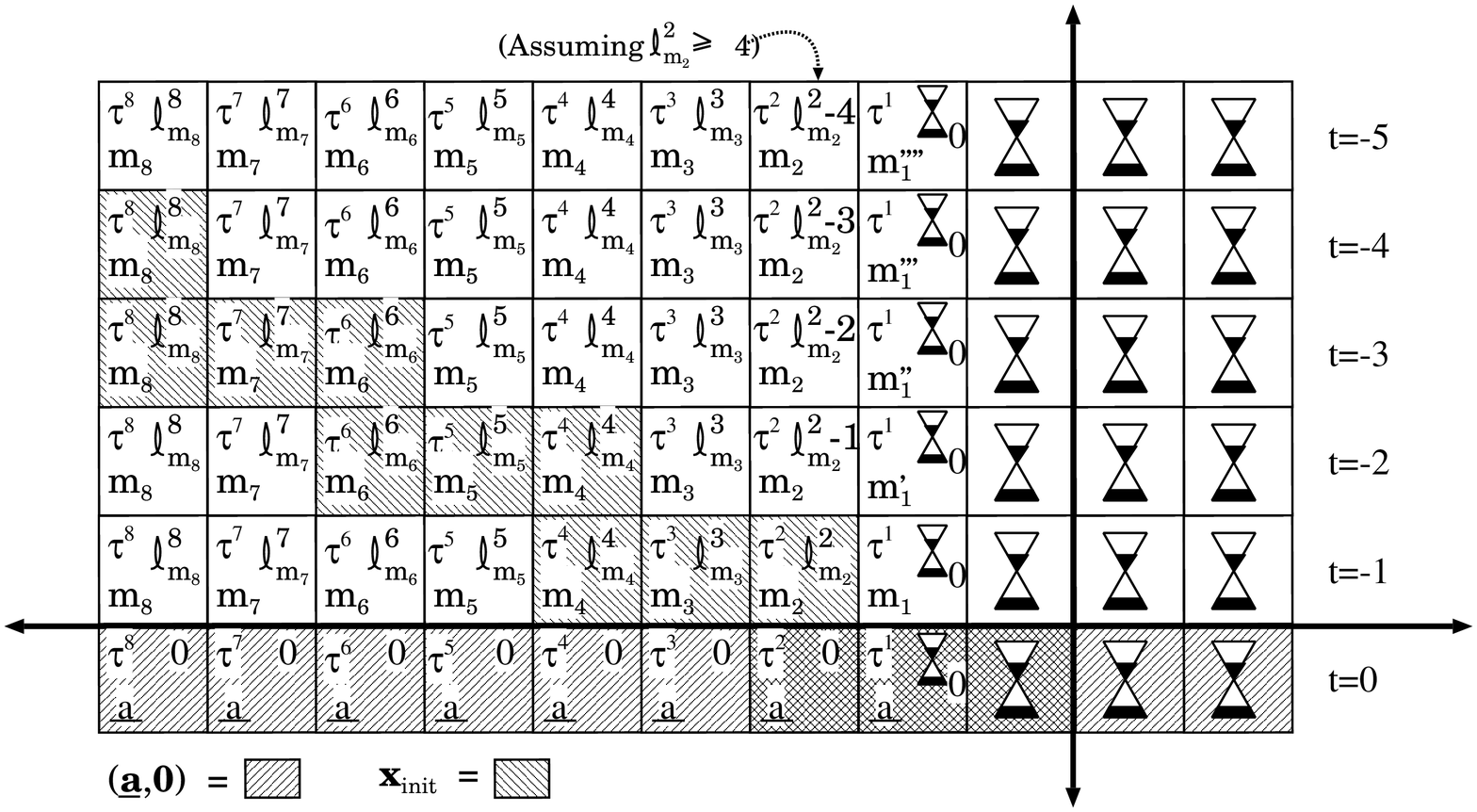}}
\caption{{\footnotesize \label{spacetimediagram}  The spacetime diagram of
 $({\bf \underline{a}},{\bf 0})$ and its five predecessors in the
adic system.  This $\mathbb{Z}^2$-indexed configuration
is an admissible element of the subshift of finite type defined using
the $2\times 2$ tiles in Figure \ref{allowedtiles}.  Finally, the tinted
diagonal ray is $\xinit$;  if read `diagonally' from top right to bottom left,
this configuration is a the spacetime diagram of $\xinit$ under the
CA $\Phi$.} }
\end{figure}

\begin{theorem}
\label{thisisit}
$(\Omega, \Phi)$ is topologically conjugate to $(X_{\cal B},V_{\cal B})$.
\end{theorem}
{\bf Proof}: 
Let $C_{1}: \Omega \rightarrow {\cal R}$ be the map which takes a
sequence $\{S_{k}\}_{k\,\in {\mathbb Z}}$ to the ray with steps
$S_{k}.$ Note that if $C_{1}({\bf x})$ is a fragment from an element
${\bf y}$ in $Y$, then $C_{1}(\Phi({\bf x}))$ is also a fragment from
$Y$, by Lemma \ref{calemma} and also by the definition of $\Phi$.   Note also (see Figure \ref{spacetimediagram})
that $C_{1}(\xinit)$ is a fragment from the
spacetime diagram of the minimal element $({\bf \underline{a}},{\bf 0})$ where
\[({\bf \underline{a}},{\bf 0}):= \ldots {\underline{a}} \stackrel{0}{\rightarrow}
{\underline{a}}\stackrel{0}{\rightarrow}\ldots 
{\underline{a}}\stackrel{0}{\rightarrow}
{\underline{a}}\stackrel{\clock_{0}}{\rightarrow}
\cdot \clock \clock \clock \ldots\, . \]

 Now define $C:=F\circ C_{1}$. By the above remarks, $C$ maps $\Omega$ into $X_{\cal B}$. Note that $C(\xinit)= 
({\bf \underline{a}},{\bf 0})    $, and that by induction
\[C(\Phi^{n}(\xinit))= V_{\cal B}^{n} (C(\xinit)),\]
for all $n\, \in {\mathbb N}$. It can be seen that $C: \{\Phi^{n}(\xinit)\}_{n\geq 0}\rightarrow 
\{V_{\cal B}^{n}(({\bf \underline{a}},{\bf 0}
))\}_{n\geq 0}$ is a uniformly continuous bijection with a uniformly continuous inverse. Thus $C$ can be extended to a 
conjugacy between $(\Omega, \Phi)$ and $(X_{\cal B},V_{\cal B})$.
\endproof

\section{\label{examples}Examples}

\subsection{Substitutions}
A {\em substitution} is a map $\btau: {\cal A} \rightarrow {\cal
A}^{+}$; we will write ${\btau}(a) = {\bf a}= a_{0}a_{2} \ldots
a_{l_a}$. We extend $\btau$ to a map $\btau: {\cal
A}^{+}\rightarrow {\cal A}^{+}$ by concatenation: if ${\bf a}=
a_{0}a_{1}\ldots a_{k}$, then $\btau({\bf a}):= \btau(a_{0})\btau(a_{1})
\ldots \btau(a_{k})$. In this way the composition $\btau\circ \btau:{\cal
A}^{+}\rightarrow {\cal A}^{+}$ is well defined; 
we will write $\btau^{\circ n}$ to refer to the
 $n$-fold composition of $\btau$.
The substitution $\btau$ is extended to a map $\btau: {\cal A}^{\mathbb
Z} \rightarrow {\cal A}^{\mathbb Z} $ defined by $\btau (\ldots
a_{-1}\cdot a_{0}\, a_{1} \ldots ) = (\ldots \btau(a_{-1})\cdot
\btau(a_{0}) \btau(a_{1}) \ldots )$. 
We say $\btau$ is {\em proper} if there exist $\underline{a},
\,\overline{a}$ in ${\cal A}$ such that for each $a\,\in {\cal A}$,
$\btau(a)$ starts with $\underline{a}$ and ends with $\overline{a}$.
 We say ${\btau}$ is
{\em primitive} if there exists a positive integer $k$ such that for
any $a\, \in {\cal A}$, all letters of ${\cal A}$ appear in
$\btau^{\circ k}(a)$ (this requires that for some letter $a$, $l_a \geq
1$). 
A {\em language} is a collection of words from 
${\cal A}^{+}$.
 If ${\bf a}\, \in {\cal A}^{\mathbb L}$,
the word ${\bf W}$ is a {\em factor of ${\bf a}$} if there is some
$n\, \in {\mathbb L}$ such that $a_{n} a_{n+1}\ldots a_{n+|{\bf W}|} =
{\bf W}$.
If ${\cal L}_{\btau}$ is the
language generated by the words $\btau^{\circ n}(a)$ with $n$ positive and $a
\, \in {\cal A}$, then let $X_{\btau}$ be the subshift associated with ${\cal
L}_{\btau}$.
  We call $X_{\btau}$  the {\em
(substitution) subshift defined by} ${\btau}$.  Henceforth we assume
that $\btau$ is primitive. In that case, $(X_{\btau}, \sigma)$ is {\em
minimal} (every $x \, \in X_{\btau}$ has a dense orbit).  For details
see \cite{pf}.  A {\em fixed point} of $\btau$ is a sequence ${\bf a}
\, \in {\cal A}^{Z}$ such that $\btau({\bf a}) = {\bf a}$. If ${\bf
a}$ is a fixed point for $\btau$, then $\btau(a _{0})$ starts with
$a_{0}$ and $\btau(a_{-1})$ ends with $a_{-1}$. Conversely, if there
exist letters $r$, $l$ such that $\btau(l)$ ends with $l$ and $\btau(r)$
starts with $r$, then ${\bf a}:=\lim_{n\rightarrow \infty}
\btau^{\circ n}(l)\cdot \btau^{\circ n}(r)$ is the unique fixed point satisfying
$a_{-1}=l$ and $a_{0}=r$. The fixed point ${\bf a}$ is {\em admissible
for $\btau$} if the word $a_{-1}a_{0}$ occurs in $\btau^{\circ n}(a)$, for
some $a\,\in {\cal A} $ and some positive $n$.  If ${\bf a}$ is
admissible, and $\btau$ is primitive, then ${\cal O}_{\sigma}({\bf
a}):=\overline{\{\sigma^{n}({\bf a})\}} = X_{\btau}$.  Using the
pigeonhole principle, there exists $n\geq 1$ such that $\btau^{\circ n}$ has
at least one admissible fixed point, and since $\btau$ and  $\btau^{\circ n}$
define the same subshift, we can assume that any primitive
substitution subshift is the orbit closure of some admissible fixed
point. If $\btau$ is proper, then it has a unique (admissible) fixed
point. The substitution $\btau$ is called {\em aperiodic} if $ {\cal
O}_{\sigma}({\bf a})$ is infinite, i.e. if ${\bf a}$ is not
$\sigma$-periodic.

\subsubsection{The Bratteli diagram associated with a proper substitution}
\begin{figure}[h]
\centerline{\includegraphics[scale=1.0]{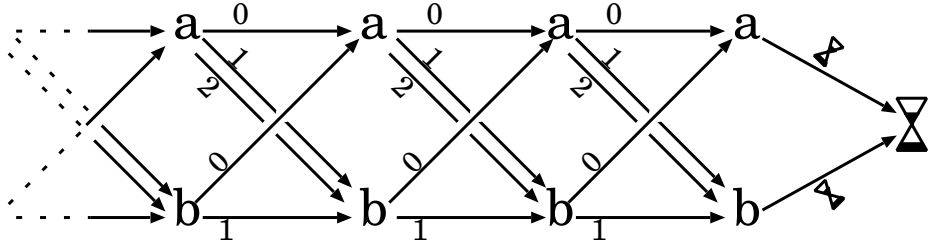}}
\caption{{\footnotesize \label{Brattelidiagram}The Bratteli diagram associated with the substitution
$\btau(a) =abb,\,\,\, \btau(b)=ab\,$.}}
\end{figure}

Let $\btau$ be a proper substitution on ${\cal A}$, with all words ${\btau}(a)$ starting with $\underline{a}$.
The {\em Bratteli diagram
associated with $\btau$} has vertex sets ${\cal V}_{n} = \{(a,n):a\, \in
{\cal A}\},$  
for each $n\geq 1$  (in the notation of Section \ref{proof}, 
we have ${\btau}^n:= \btau$ for all $n\in\mathbb{N}$).
   There is exactly one edge from each vertex in ${\cal V}_{1}$ to $\clock$.
If $1 \leq m<n$, then the
number of paths in ${\cal B}$ from $(a,n) $ to $(b,m)$ is the number
of occurences of $b$ in ${\btau}^{\circ (n-m)}(a).$ If $\btau$ is primitive,
then there is a positive $k$ such that for any two letters $a$ and
$b$, there is at least one path from $(a,n+k)$ to $(b,n)$. The
Bratteli diagram ${\cal B}' =({\cal V}',{\cal E}')$ for $\btau^{\circ k}$ is
the telescoping of the Bratteli diagram ${\cal B}= ({\cal V},{\cal
E})$ for $\btau$, with ${\cal V}_{n}'= {\cal V}_{nk}$ for $n>1$ and
${\cal V}_{1}' = {\cal V}_{1}$. Thus, Bratteli diagrams associated with proper
 substitutions are simple,  focused and of bounded width. An example of such a Bratteli diagram is given in Figure \ref{Brattelidiagram}.

\subsection{Odometers}

Let ${\bf q} : = (q_{1}, q_{2}, \ldots )$ be an ordered set of integers $ \geq 2$ (the {\em quotient
set}). Let ${\cal Z}({\bf q}):= \prod_{l = 1}^{\infty} {\mathbb
Z}_{q_{l}}$ be the Cartesian product set. Then  $({\cal Z}({\bf q}),
\oplus)$ is a (compact, zero-dimensional, profinite, abelian) group where ``$\oplus$'' is defined
as addition ``with carry'': if ${\bf x}=( \ldots , x_{2},x_{1})$ and 
${\bf y}=( \ldots ,y_{2}, y_{1})$, then ${\bf
x}\oplus {\bf y} := ( \ldots,r_{2} ,r_{1})$ where
$x_{1}+y_{1} = k_{1}q_{1}+r_{1}$ for some $k_{1}\geq 0$ and $0 \leq r_{1}<q_{1}$,  and 
for each $n\geq 2$,  
 \[ k_{n-1}+x_{n}+y_{n} = k_{n}q_{n}+r_{n}, \]
with $k_{n}\geq 0$ and $0\leq r_{n} < q_{n}$ . 
 If $p$ is prime, then the {\em multiplicity of p in} ${\bf q}$ is
the sum number of times (possibly infinite) that $p$ occurs in the
prime decomposition of all the elements of sequence ${\bf q}$. 

Let 
${\bf 1}:= ( \ldots ,0,0,1)$. We define the
{\em odometer} $\varsigma: {\cal Z}({\bf q})
\rightarrow {\cal Z}({\bf q})$ as $\varsigma ({\bf g}) =
{\bf g}\oplus {\bf 1}$. In other words, an odometer is a group
rotation on ${\cal Z}({\bf q})$. 

\begin{theorem}
\label{downarowicz}
$({\cal Z}({\bf q}), \varsigma) $ and $({\cal Z}({\cal
Q^{*}}), \varsigma)$ are topologically conjugate if and only
if  every prime $p$ has equal multiplicity in ${\bf q}$ and
${\cal Q^{*}}$.
\hfill {\rm\cite[Thm 1.2]{down}}
\end{theorem}

Thus we can assume that elements in ${\bf q}$ are prime. 
 If $q_{i}=q$ for each $i$ then we will write ${\cal Z}( q)$ for
${\cal Z}({\bf q})$; this is know as the ``$q$-adic'' odometer, and is a model for ``base $q$'' arithmetic. Odometers have an adic representation where each ${\cal V}_{n}$ is a one-element set and for each $n$, there are $q_{n+1}$ edges in ${\cal E}_{n+1}$. Adic systems are often seen as ``generalised odometers''.

The connection between adic systems, substitutions and many odometers is given by the next result:

\begin{theorem}
\label{sub=bratt}
Let $\btau $ be a primitive substitution, with fixed point ${\bf x}.$
\begin{enumerate}
\item If $\btau$ is aperiodic, then $(X_{\btau},\sigma, {\bf x})$ is
pointedly isomorphic to $(X_{\cal B},V_{\cal B}, {\bf x}_{\min})$, 
for some Bratteli diagram arising from a proper aperiodic substitution.
\item If $\btau$ is periodic, then $(X_{\cal B},V_{\cal B})$ is
isomorphic to the odometer $({\cal Z}(\ldots,k,k,p), \varsigma)$ where
$p$ is the periodicity of $X_{\btau}$, and $k\, \in {\mathbb N}$.
\end{enumerate}
\hfill\rm{
\cite[Thm 17]{f} or \cite[Prop 20]{dhs}}
\end{theorem}

\begin{corollary}
\label{cor1}
\begin{enumerate}

\item
Aperiodic, primitive substitution systems can be embedded in a CA.
\item
Odometers whose quotient set has finitely many primes can be embedded in a CA.

\end{enumerate}

\end{corollary}

{\bf Proof:} The first statement follows from Theorem \ref{sub=bratt}.

To prove the second statement: Suppose that
${\bf q}_{f}$ and  ${\bf q}_{i}$ are the set of  primes
in ${\bf q}$ with finite and infinite  multiplicity respectively.
Let $N:=(\prod_{p\,\in {\bf q}_{f}} p)$ and $M:=(\prod_{p\,\in {\bf q}_{i}} p)$. Theorem \ref{downarowicz} tells us that we can assume that ${\bf q} = (\ldots, M,M, N)$.
If $\btau$ is the substitution on a $1$-letter alphabet ${\cal
A}=\{a\}$ defined by $\btau(a) = a^{M}$, then the Bratteli system $(X_{\cal B},V_{\cal
B})$ associated to $\btau$ is topologically isomorphic to $({\cal
Z}(M),\varsigma)$, by Theorem \ref{sub=bratt}. If the Bratteli diagram for $\btau$ is modified by
putting $N$ edges from the vertex $a$ in  ${\cal V}_{1}$ to $\clock$,
then the resulting Bratteli system is topologically isomorphic to
$({\cal Z}(\ldots M, M, N),\varsigma)$.
\endproof

\subsection{Toeplitz subshifts}
If ${\cal A}$ is a finite alphabet, and $X= {\cal A}^{\mathbb Z}$, a
{\em Toeplitz sequence} is a $\sigma$-aperiodic element ${\bf x} \,
\in X$ such that for each $m \, \in \mathbb Z$, there is some $p$ such
that $x_{m+pn}=x_{m}$ for each $n \, \in \mathbb Z$. The {\em Toeplitz
subshift} associated with ${\bf x}$ (also sometimes called the {\em Toeplitz flow}) is the subshift $({\cal
O}_{\sigma}({\bf a}), \sigma ).$ See \cite{down}.

Given a Bratteli diagram ${\cal B}$ whose vertices are labelled from
${\cal A}$, we say that ${\cal B}$ has the {\em equal path number
property} if for each $n \geq 1$, $l^n_a$ is independent of
$a$. A dynamical system $(Y,T)$ with $Y$ a metric space  is {\em expansive} if there exists
some $\epsilon >0$ such that for each ${\bf x} \neq {\bf y}$, there is some integer $n$ such that $d(T^{n}({\bf x}), T^{n}({\bf y}))\geq \epsilon$.

\begin{theorem}
\label{toeplitzthm}
The family of expansive adic systems associated to properly ordered Bratteli diagrams with the equal path number property
is upto conjugacy, the family of Toeplitz flows.
\hfill {\rm\cite[Thm8]{gj1}}
\end{theorem}

\begin{lemma}
\label{focusedtoeplitzlemma}
Given a Toeplitz system $(Y, \sigma)$, the Bratteli diagram constructed in Theorem \ref{toeplitzthm} is focused.
\end{lemma}

{\bf Proof:} The proof of Theorem \ref{toeplitzthm} involves finding a
sequence of collections of words $\{ {\cal W}_{n} \}_{n}$, and a
sequence of words ${\bf W}_{n} \, \in {\cal W}_{n}$ such that all
words in ${\cal W}_{n+1}$ are concatenations of words from ${\cal
W}_{n}$, and also begin with the word ${\bf W}_{n}$. In fact if ${\bf
x}$ is the Toeplitz sequence, there will be a sequence $p_{n}$ such
that ${\bf W}_{n} = x_{0}x_{1} \ldots x_{p_{n}-1}$. The constructed
Bratteli diagram will have vertices in ${\cal V}_{n}$ labelled by
words in ${\cal W}_{n}$, and, if ${\bf W}\, \in {\cal W}_{n+1}$ is
such that it is a concatenation ${\bf W} = {\bf W}^{1}{\bf W}^{2}
\ldots {\bf W}^{k}$ of words in ${\cal W}_{n}$, there will be an
ordered edge with label $i$ from the vertex representing ${\bf W}$ in
${\cal V}_{n+1}$ to the vertex in ${\cal V}_{n}$ representing ${\bf
W}^{i}$.  Since all words in ${\cal W}_{n+1}$ begin with ${\bf
W}_{n}$, this means that all minimal edges in ${\cal E}_{n+1}$ have
the same range.
\endproof

One can construct a Toeplitz sequence ${\bf x}$ in the following way,
which is described in \cite{down}. Fix an alphabet ${\cal A}$.  Pick
an integer $s_{1}\geq 2$, and pick a subset $F_{1} \subset [0,s_{1})$
and a ${\cal A}$-valued labelling of $F_{1}$, say $A_{1}\, \in {\cal
A}^{F_{1}}$. For each $n \, \in \mathbb Z$, let $x_{ns_{1} +F_{1}} =
A_{1}$. Let ${\rm Per}_{s_1}({\bf x}) :=F_{1}+s_{1}{\mathbb Z}$.  Given
$s_{k}$, choose $s_{k+1}$ such that $s_{k+1}$ is an essential multiple
of $s_{k}$. Choose $F_{k+1} \subset [0,s_{k+1}-1)\backslash
{\rm Per}_{s_{k}}({\bf x})$, and choose $A_{k+1} \, \in {\cal
A}^{F_{k+1}}$. Let $x_{ns_{k+1} +F_{K+1}} = A_{k+1}$ for all integers
$n$. For any integer $n$, the sets $[ns_{k},(n+1)s_{k})$ are {\em
k-intervals}, and the co-ordinates in a $k$-interval not filled in
during the first $k$ steps in the construction are called {\em k holes}.
The choices above must be made so that eventually all of ${\bf x}$ is
filled in. Also, after completing this procedure, we redo the
construction, at each stage $k$ re-defining $F_{k}$ so that all
essentially $s_{k}$-periodic parts of ${\bf x}$ are filled in.

The next lemma identifies which properly ordered Bratteli diagrams with the equal path number property have bounded width:

\begin{lemma}
\label{boundedwidthlemma}
Let ${\bf x}$ be a Toeplitz sequence constructed in the above
manner. Suppose that for each $k$, $s_{k+1}/ s_{k}\leq K$, and also so
that there are at most $log_{|{\cal A}|}(K)$ $k$-holes. Then the
Bratteli diagram constructed in the proof of Theorem \ref{toeplitzthm}
is of bounded width $K$.
\end{lemma}
 {\bf Proof:} If there are $log_{|\cal A|}(K)$ $k$-holes, this means
 that there are at most $K$ elements in ${\cal W}_{k}$, where the sets
 ${\cal W}_{k}$ are those defined in the proof of Proposition
 \ref{focusedtoeplitzlemma}. This bounds the number of vertices in
 ${\cal V}_{k}$. If $s_{k+1}/s_{k} \leq K$, this means that each word in
 ${\cal W}_{k+1}$ are concatenations of at most $K$ words from ${\cal W}_{k}$. Thus there at most $K$ edges emanating from each vertex in ${\cal W}_{k+1}$.
\endproof
\begin{corollary}
\label{toeplitzcorollary}
Toeplitz systems satisfying the conditions of Lemma \ref{boundedwidthlemma}
can be embedded in a cellular automaton.
\end{corollary}

{\bf Remark}: Some non-primitive substitutions systems also have
adic representations, which are focused and of bounded
width. For example, the {\em Chacon substitution}
$\btau(a)=aaba,\,\,\, \btau(b)=b$ is conjugate to an {\em induced
system} of the substitution $\btau^*(a)=aab, \,\,\, \btau^*(b)=abb$
(\cite[ \S4.2]{gj1}). Which aperiodic substitutions have a focused bounded width
representation? The same can be asked of  {\em finite rank } transformations obtained
using the ``cutting and stacking'' method with a bounded number of
cuts and spacers: do they have a focused bounded width Bratteli representation?

{\footnotesize
\bibliographystyle{alpha}
\bibliography{bibliography}

\begin{thebibliography}{CPY07}

\bibitem[CPY07]{cpy}
Ethan~M. Coven, Marcus Pivato, and Reem Yassawi.
\newblock Prevalence of odometers in cellular automata.
\newblock {\em Proc. Amer. Math. Soc.}, 135(3):815--821 (electronic), 2007.

\bibitem[CY07]{cy}
Ethan Coven and Reem Yassawi.
\newblock Embedding odometers in cellular automata.
\newblock {\em (preprint)}, 2007.

\bibitem[DHS99]{dhs}
F.~Durand, B.~Host, and C.~Skau.
\newblock Substitutional dynamical systems, {B}ratteli diagrams and dimension
  groups.
\newblock {\em Ergodic Theory Dynam. Systems}, 19(4):953--993, 1999.

\bibitem[Dow05]{down}
Tomasz Downarowicz.
\newblock Survey of odometers and {T}oeplitz flows.
\newblock In {\em Algebraic and topological dynamics}, volume 385 of {\em
  Contemp. Math.}, pages 7--37. Amer. Math. Soc., Providence, RI, 2005.

\bibitem[Fog02]{pf}
N.~Pytheas Fogg.
\newblock {\em Substitutions in dynamics, arithmetics and combinatorics},
  volume 1794 of {\em Lecture Notes in Mathematics}.
\newblock Springer-Verlag, Berlin, 2002.
\newblock Edited by V.\ Berth\'e, S.\ Ferenczi, C.\ Mauduit and A.\ Siegel.

\bibitem[For97]{f}
A.~H. Forrest.
\newblock {$K$}-groups associated with substitution minimal systems.
\newblock {\em Israel J. Math.}, 98:101--139, 1997.

\bibitem[GJ00]{gj1}
Richard Gjerde and {\O}rjan Johansen.
\newblock Bratteli-{V}ershik models for {C}antor minimal systems: applications
  to {T}oeplitz flows.
\newblock {\em Ergodic Theory Dynam. Systems}, 20(6):1687--1710, 2000.

\bibitem[GM05]{gm}
William Geller and Micha{\l} Misiurewicz.
\newblock Irrational life.
\newblock {\em Experiment. Math.}, 14(3):271--275, 2005.

\bibitem[Hed69]{h}
G.~A. Hedlund.
\newblock Endormorphisms and automorphisms of the shift dynamical system.
\newblock {\em Math. Systems Theory}, 3:320--375, 1969.

\bibitem[HPS92]{hps}
Richard~H. Herman, Ian~F. Putnam, and Christian~F. Skau.
\newblock Ordered {B}ratteli diagrams, dimension groups and topological
  dynamics.
\newblock {\em Internat. J. Math.}, 3(6):827--864, 1992.

\end{thebibliography}
}

\end{document}